\documentclass[preprint]{elsarticle}

\usepackage[english,francais]{babel}

\usepackage{amsmath, amssymb, amsfonts, amstext}
\usepackage{multicol}
\usepackage{fancyhdr}
\usepackage{natbib}
\usepackage[latin1]{inputenc}
\usepackage[T1]{fontenc}
\usepackage{latexsym}
\usepackage{graphicx}
\usepackage{stmaryrd}        
\usepackage{dsfont}          
\usepackage{subfigure, pgf}  
\usepackage[center]{caption} 
\usepackage{mathrsfs}        
\usepackage{amsthm}          
\usepackage{float}           
\usepackage{booktabs}        
\usepackage[normalem]{ulem}  
\usepackage{todonotes}       
\usepackage{verbatim} 
\usepackage{cancel}          
\usepackage{cases}           
\usepackage{multicol}        
\usepackage{setspace}        

\newtheorem{theorem}{Theorem}[section]

\newtheorem{e-proposition}[theorem]{Proposition}

\newtheorem{e-definition}[theorem]{Definition\rm}

\newcommand{\ds}{\displaystyle}

\newcommand{\R}{\mathbb{R}}

\newcommand{\bx}{{\mathbf{x}}}
\newcommand{\by}{{\mathbf{y}}}

\newcommand{\sA}{\mathsf{A}}

\newcommand{\sC}{\mathsf{C}}

\newcommand{\epss}{\varepsilon}
\newcommand{\promille}{\small \hbox{$\,^0\!/_{00}$}}

\newenvironment{systeme*}{\begin{equation*}
\left\{ \begin{aligned}}{\end{aligned} \right.
\end{equation*}}

\journal{Applied Mathematics Letters}

\begin{document}

\begin{frontmatter}




\selectlanguage{english}
\title{Explicit computation of the electrostatic energy for an elliptical charged disc}


\author[imt,cerfacs]{S. Laurens}
\ead{sophie.laurens@insa-toulouse.fr}
\address[imt]{Mathematical Institute of Toulouse, 118 route de Narbonne, F-31400 Toulouse}
\address[cerfacs]{CERFACS, 42 Avenue Gaspard Coriolis F-31100 Toulouse}
\author[inria]{S. Tordeux}
\address[inria]{INRIA Bordeaux Sud-Ouest-LMA, avenue de l'Universit\'e, F-64013 Pau}

\begin{abstract}
This letter describes a method for obtaining an explicit expression for the electrostatic energy of a charged elliptical infinitely thin disc. The charge distribution is assumed to be polynomial. Such explicit values for this energy are fundamental for assessing the accuracy of boundary element codes. The main tools used are an extension of  Copson's method and a diagonalization, given by Leppington and Levine, of the single-layer potential operator associated with the electrostatic potential created by a distribution of charges on the elliptical disc. 
\end{abstract}



\end{frontmatter}


\selectlanguage{english}


\vspace{-1.5em}

\section{Introduction}

 In recent years, integral equations have become an essential tool for solving both industrial and scientific problems in  electromagnetism and acoustics. The assessment of the accuracy delivered by such codes, in particular in their handling of the  singular integrals involved, is a  major issue.  Here, we present a  method for deriving an analytical expression for the electrostatic energy of a charged elliptical infinitely thin plate, providing a means for the validation of these codes. \\

Let us denote by $\ds   \sA = \left\{ (x_1,x_2) \in \R^ 2 \mbox{ with }  x_1^{2}/a^{2}+x_2^{2}/b^{2}-1<0  \mbox{ and  } a>b\right\}$ the ellipse with  major and minor semi-axes $a$ and $b$. Let $f$  be the electrostatic potential generated by a density of charges $\sigma$ distributed over $\sA$: 
 \begin{equation} \label{eq:1}
 f(\bx) = \dfrac{1}{4\pi}\int_{\sA}\dfrac{\sigma(\by)}{\vert \bx-\by \vert}ds_\by  \qquad \mbox{ for all } \bx \mbox{ in }  \sA,  \vspace{-0.5em}
\end{equation}
with $\bx = (x_1,x_2)$. The units have been chosen such that the electric permittivity of air is $1$. The electrostatic energy $I$ can be expressed in either the two following forms:
\begin{equation}\label{eq:2}
  I_{\sigma} = \int_{\sA}  f(\bx) \: \overline{\sigma(\bx)} ds_\bx =   \int_{\sA} \dfrac{1}{4\pi}\int_{\sA}\dfrac{\overline{\sigma(\bx)}\:\sigma(\by)}{\vert \bx-\by \vert} \: ds_\bx ds_\by .\vspace{-0.5em}
\end{equation}
We aim in this letter at proving  and  numerically illustrating the following theorem, where $\epss$ is the eccentricity of the ellipse $\sA$ given by $\epss =\ds \sqrt{1 - b^2/a^2}$.
\begin{theorem}\label{theo}
Let $\sigma(\bx) = \alpha_0+\alpha_1 x_1/a+ \alpha_2 x_2/b$, with $\alpha\in \R^3$, be the distribution of charges over $\sA$. The corresponding electrostatic energy is given by
\begin{equation*}
    I_{\sigma}=  \dfrac{8ab^{2}}{15 \pi }   \left[\left(5\alpha_0^2+ \alpha_2^2  \right) K(\epss) + \left(\alpha_1^2-\alpha_2^2 \right) \dfrac{K(\epss)- E(\epss)}{\epss^2}\right],
\end{equation*}
with $K(\varepsilon)$ and $E(\epss)$ the complete  elliptic integrals of the first and second kind  
  \begin{equation}
  \label{eq:KEeps}\ds K(\epss) = \int_0^{\pi/2} \dfrac{d\phi}{\sqrt{ 1-\epss^2 \sin^2 \phi }}\quad\hbox{ and }\quad\ds E(\epss) = \int_0^{\pi/2} \sqrt{ 1-\epss^2 \sin^2 \phi} \: d\phi.
  \end{equation}
\end{theorem}

\section{Diagonalization of the electrostatic energy}

Following \cite{leppington1973}, we consider the spheroidal coordinate system  $(\theta, \varphi)$ giving a parametrization of $\sA$ in terms of the unit half-sphere 
\begin{equation}\label{eq:3}
  x_1=a\, \sin\theta \, \cos\varphi\: \mbox{ and } \: x_2=b\, \sin\theta\,  \sin\varphi,  \mbox{ with  }\:  \theta \in [0,\pi/2], \: \varphi \in [0,2\pi].
\end{equation}
The elemental area associated with the new variables is $ ab  \cos\theta \sin\theta  \,d\theta d \varphi$. In these spheroidal coordinates, the electrostatic potential $f$ defined in \eqref{eq:1} and the electrostatic energy can also be written in terms of $\theta$ and $\varphi$ as
\begin{equation*}
\left\{
\begin{array}{llll}
   \ds f(\theta,\varphi) = \dfrac{ab}{4\pi}\int_{0}^{\pi/2}\int_{0}^{2\pi}\dfrac{g(\theta',\varphi')}{d(\theta,\varphi,\theta',\varphi')}  \sin \theta' \: d\theta'd\varphi',\\
    \ds I_{\sigma}
 =   ab\int_{0}^{\pi/2}  \int_{0}^{2\pi} \overline{f(\theta,\varphi) } g(\theta,\varphi)  \sin \theta  \: d\theta d\varphi,
\end{array}
\right.
\end{equation*}
with $\ds g(\theta,\varphi) = \sigma(\theta,\varphi) \; \cos\theta$ and $d(\theta,\varphi,\theta',\varphi')$ the distance separating $\bx$ from $\by$ $d = \vert \bx-\by \vert $ expressed in the spheroidal coordinates \eqref{eq:3}. \\
The next step consists  in introducing a well chosen spectral basis 
for the half-sphere involving the even Legendre functions $Q^{m}_{n}$ normalized by
\begin{equation}
\int_{0}^{\pi} Q_{n}^{m} \left(  \cos\theta\right)  Q_{n'}^{m}\left(  \cos\theta\right)  \sin \theta d\theta= \delta_{n,n'}.
\end{equation}
This basis yields a block diagonalization of  the convolution operator (see \cite{leppington1973})
\begin{multline}
  \label{eq:8}
\frac{1}{d}=\frac{1}{\sqrt{ab}} \:\sum_{n=0}^{\infty}\sum_{\substack{ m = - n \\ n-m \text{ even}}}^{n}
\sum_{\substack{ m' = - n \\ n-m' \text{ even}}}^{n}d_{mm^{\prime}}^{n}Q_{n}^{m}\left(  \cos\theta\right)
Q_{n}^{m^{\prime}}\left(  \cos\theta^{\prime}\right)  e^{i\left(
m\varphi-m^{\prime}\varphi^{\prime}\right)  },\\ 
 \text{with }
  d_{mm^{\prime}}^{n}=\frac{Q_{n}^{m}\left(  0\right)  Q_{n}^{m^{\prime}}\left(
0\right)  }{2n+1}\int_{0}^{2\pi}\frac{e^{i(m-m^{\prime})\varphi}}{\sqrt{\frac
{b}{a}\cos^{2}\varphi+\frac{a}{b}\sin^{2}\varphi}}d\varphi.
\end{multline}
The functions $f$ and $g$ can be expanded in this basis as
\begin{equation}\label{eq:10}
   u(\theta,\varphi)=\ds\sum_{n=0}^{\infty} \: \sum_{\substack{ m = - n \\ n-m \text{ even}}}^n u_{n}^{m}Q_{n}^{m}\left(  \cos\theta\right) e^{im\varphi}\quad\hbox{ with }u=f\hbox{ or }g
\end{equation}
with
\begin{equation}\label{eq:7bis}
 u^{m}_{n}\;=\;\frac{1}{\pi}\ds\int_0^{\pi/2}\int_0^{2\pi} u(\theta,\varphi)\;Q_n^m(\cos\theta)\;e^{-im\varphi}\sin(\theta)d\theta d\varphi.
\end{equation}
Due to \eqref{eq:8}, coefficients $f^m_{n}$ are related to $g^m_{n}$ by
\begin{equation}
f^{m}_{n}\;=\;\dfrac{\sqrt{ab}}{4} \sum_{\substack{ m' = - n \\ n-m' \text{ even}}}^{n} d_{mm^{\prime}}^{n} g_{n'}^{m^{\prime}}.
\end{equation}
Moreover, the orthogonal properties of the spectral basis yield 
\begin{equation*}
I_{\sigma}\;=\; \dfrac{\pi}{4}(ab)^{3/2}  \sum_{n,m,m'} d_{mm^{\prime}}^{n}  g_{n}^{m} \overline{g_{n}^{m^{\prime}}},
\end{equation*}
where we have lightened the notation by making  the range of the summation index implicit. Indices $n, n'$ are varying from $0$ to $\infty$, and $m,m'$ are such that $\vert m \vert \leq n$, $\vert m' \vert \leq n'$, with $n-m$ and $n-m'$ even.
Substituting expression   \eqref{eq:8} for $d_{mm^{\prime}}^{n} $ and introducing the eccentricity of the ellipse $\epss$, we get
\begin{equation}\label{eq:I}
    I_{\sigma} =  \pi\:  ab^{2} \:  \sum_{n,m,m'} g_{n}^{m} \: \overline{g_{n}^{m^{\prime}}}\: \: \dfrac{Q_{n}^{m}\left(  0\right) Q_{n}^{m^{\prime}}\left(0\right)  }{2n+1} \: \int_{0}^{\pi/2}\dfrac{\cos( m - m')\varphi}{\sqrt{ 1-\epss^2 \cos^2 \varphi }}d\varphi.
\end{equation}

\section{Electrostatic energy for an affine distribution of  charges}
This section is dedicated to the calculation of the electrostatic energy gene\-rated by an affine density of charges.\\\noindent
\textbf{Proof of Theorem \ref{theo}.}
Let $\sigma_0(\bx)=1$, $\sigma_1(\bx)=x_1/a$ and $\sigma_2(\bx)=x_2/b$. In view of the symmetry of $\sA$ with respect to $x_1$ and $x_2$, we have
\begin{equation*}
 \int_{\sA} \dfrac{1}{4\pi}\int_{\sA}\dfrac{\sigma_i(\bx)\sigma_j(\by)}{\vert \bx-\by \vert} \: ds_\bx ds_\by\;=\;0\quad \quad\hbox{ for }i\neq j.
\end{equation*}
Consequently, the electrostatic energy $I_{\sigma}$ can be expanded as
\begin{equation*}
 I_{\sigma} =     \alpha_0^2 \: I_{\sigma_0} +\alpha_1^2 \: I_{\sigma_1} +\alpha_2^2 \: I_{\sigma_2} 
\end{equation*}
The result will follow from the computation of $I_{\sigma_0}$, $I_{\sigma_1}$ and $I_{\sigma_2}$.

  \subsection{Computation of $I_{\sigma_0}$}

For $\sigma(\bx) = \sigma_0(\bx) = 1$, the function $\ds g(\theta, \varphi) = \cos \theta$ does not depend on $\varphi$. The $g_n^m$ coefficients  are independent from $a$ and $b$, and since $g_n^m=0$ for all $m\neq0$, the function $g$ can be expanded as $g(\theta)  = \ds \sum_{n=0}^{+\infty} g_{n}^0 Q_{n}^{0}\left(\cos\theta\right)$ .
Due to \eqref{eq:I}, the electrostatic energy $ I_{\sigma_0}$ depends only on $a$ and $b$ and is given by  
\begin{equation}\label{eq:12bis}
  I_{\sigma_0}(a,b) =  \pi \: ab^{2}  \: \sum_{n=0}^{+\infty}   \dfrac{\vert g_{n}^{0} \vert ^2 \: \left(Q_{n}^{0}\left(  0\right) \right)^2  }{2n+1}  \int_0^{\pi/2} \dfrac{d\varphi}{\sqrt{ 1-\epss^2 \cos^2 \varphi }}  = \kappa \: ab^{2} \: K(\epss),
\end{equation}
with $\varepsilon=\sqrt{1-b^2/a^2}$ and $\kappa$ a constant depending neither on $a$ nor on $b$.
The constant $\kappa$ is deduced from the classical case of an unit circle
which has been detailed for example in \cite{laurens2012} 
\begin{equation}\label{eq:12ter}
 I_{\sigma_0}(1,1)\;=\;\int_{\sC} \dfrac{1}{4\pi}\int_{\sC}\dfrac{1}{\vert \bx-\by \vert} \: ds_\bx ds_\by = 4/ 3.
\end{equation}
 Comparing \eqref{eq:12bis} and \eqref{eq:12ter}, this yields to $\kappa=8/3\pi$ since $K(0)=\pi/2$. Therefore,
\begin{equation}
  \label{eq:Isigma0}
 I_{\sigma_0}   =  \dfrac{8}{3 \pi } \: ab^{2} \:  K(\epss).
\end{equation}

\subsection{Computation of $I_{\sigma_1}$}

For $\sigma(\bx) = \sigma_1(\bx)=x_1/a$, the function $\ds g$ is given by $\ds g(\theta, \varphi) = \sin \theta \cos \theta \cos \varphi $. In that case, the $g_n^m$ are zero except for $\vert m \vert =1$. By definition of the Legendre functions, we have $Q_n^{-1} =- Q_n^1$. As the function $g$ is even, it emerges that $g_n^{-1} = - g_n^1$, and thus
\begin{equation*}
    I_{\sigma_1} \! =  \pi  a b^{2}   \sum_{n}    \dfrac{\left( g_{n}^{1} Q_{n}^{ 1} (0) \right)^2  }{2n+1} \! \left[ 2 \! \int_0^{\pi/2} \hspace{-.75em}  \dfrac{d\varphi}{\sqrt{ 1-\epss^2 \cos^2 \varphi }}  + 2 \! \int_{0}^{\pi/2} \hspace{-.5em} \dfrac{\cos 2\varphi \: \: d\varphi}{\sqrt{ 1-\epss^2 \cos^2 \varphi } } \right].
\end{equation*} 
Due to \eqref{eq:KEeps}, it emerges that
\begin{equation}\label{inter:1}
  I_{\sigma_1}(a,b)= \kappa \: a b^2 \:  \dfrac{K(\epss) -  E(\epss)}{\epss^2}.
\end{equation} 
To determine the constant $\kappa$, we consider again the case of an unit circle.  In this case, $I_{\sigma_1}$ can be explicitly computed (see the \ref{app}), and is given by $I_{\sigma_1}(1,1) =  2 / 15$. Evaluating \eqref{inter:1} at $a=b=1$ we get
\begin{equation}\label{I:1}
 I_{\sigma_1}(1,1) = \kappa\dfrac{ \pi}{ 4 }\quad\hbox{ since }\lim_{\varepsilon\rightarrow0}\dfrac{K(\epss) - E(\epss)}{ \epss^2}=\frac{\pi}{4}.
\end{equation}
It follows that $\kappa=8/ 15\pi$ and therefore we have
\begin{equation}
    \label{eq:Isigma1}
    I_{\sigma_1} =  \dfrac{8}{15 \pi } \:  ab^{2}  \:   \dfrac{K(\epss) - E(\epss)}{ \epss^2}.
  \end{equation}

  \subsection{Computation of $I_{\sigma_2}$}
For $I_{\sigma_2}$, we consider $\sigma(\bx) = \sigma_2(\bx)=x_2/b$, meaning that $\ds g(\theta, \varphi) =\sin \theta \cos \theta \sin \varphi$. We still have $g_n^m=0$ except for $\vert m\vert =1$,  but in that case, $g_n^{-1} =  g_n^1$. Thus
\begin{equation*}
    I_{\sigma_2} \! =  \pi  a b^{2}   \sum_{n=0}^{+\infty}    \dfrac{\left( g_{n}^{1} Q_{n}^{ 1} (0) \right)^2  }{2n+1} \! \left[ 2 \! \int_0^{\pi/2} \hspace{-.75em}  \dfrac{d\varphi}{\sqrt{ 1-\epss^2 \cos^2 \varphi }}  - 2 \! \int_{0}^{\pi/2} \hspace{-.5em} \dfrac{\cos 2\varphi \: \: d\varphi}{\sqrt{ 1-\epss^2 \cos^2 \varphi } } \right]
\end{equation*} 
Moreover,  both integrals $ I_{\sigma_1}$ and $ I_{\sigma_2}$ are equal on $\sC$ by symmetry. We obtain
\begin{equation}\label{eq:I22}
  I_{\sigma_2}=  \dfrac{8}{15 \pi } \: a b^2 \: \left( K(\epss) -  \dfrac{K(\epss) -  E(\epss)}{\epss^2} \right).
\end{equation}

\section{Numerical tests and conclusion}

Tables  \ref{table:1} and \ref{table:2} give a comparison of the exact values given by an analytical expression with numerical approximate values obtained by the boundary element code CESC of CERFACS with $P^1$ continuous elements.
 It can be observed that the two values coincide at least up to the fourth decimal digit. Table \ref{table:3} shows the maximum relative error for each of the cases of Tables  \ref{table:1} and \ref{table:2} cases, which is less than $0.35$ per mil.

\begin{table}[H]
  \centering
  \begin{tabular}[H]{ccccccc}
$a$ & 0.5 & 0.7 & 0.9 & 1.1 & 1.3 & 1.5 \\
\hline\\[-0.5em]
$I^\textrm{comp}_{\sigma_0}$ &  0.1666 &  0.2741 & 0.3939  & 0.5234 & 0.6608 & 0.8048 \\
$I^\textrm{exact}_{\sigma_0}$ &  0.1666 &  0.2741 & 0.3939  & 0.5234 & 0.6608 & 0.8048 \\
\hline\\[-0.5em]
$I^\textrm{comp}_{\sigma_1}$  $ \times 10^{-1}$& 0.0417  & 0.1455  & 0.3651 & 0.7535  & 1.3715 & 2.2781  \\
$I^\textrm{exact}_{\sigma_1}$ $ \times 10^{-1}$ & 0.0417  & 0.1456  & 0.3656 & 0.7543  & 1.3717 & 2.2800  \\
\hline\\[-0.5em]
$I^\textrm{comp}_{\sigma_2}$ $ \times 10^{-2}$ & 0.4167 & 0.6280  & 0.8426  & 1.0585  & 1.2748  & 1.4910 \\
$I^\textrm{exact}_{\sigma_2}$ $ \times 10^{-2}$ &0.4167 & 0.6280  & 0.8427  & 1.0586  & 1.2748  & 1.4911
  \end{tabular}
  \caption{\small Exact and computed values of the electrostatic energy, $I^\textrm{exact}_{\sigma_i}$  and  $I^\textrm{comp}_{\sigma_i}$,  for an elliptical disc with minor axis $b = 0.5$, given by  \eqref{eq:Isigma0},  \eqref{eq:Isigma1},  \eqref{eq:I22} .}
  \label{table:1}
\end{table}  

\vspace{-1.em}

\begin{table}[H]
  \centering
  \begin{tabular}[H]{ccccccc}
$a$ & 0.75 & 0.9 & 1.05 & 1.2 & 1.35 & 1.5 \\
\hline\\[-0.75em]
 $I^\textrm{comp}_{\sigma}$  & 2.7734 & 3.6157 & 4.5162 &  5.4706 & 6.4761 & 7.5313 \\
 $I^\textrm{exact}_{\sigma}$ & 2.7736 & 3.6159 & 4.5165 &  5.4708 & 6.4763 & 7.5316 
  \end{tabular}
  \caption{\small Affine density of charges $\sigma=x_1+2x_2+3$  for an elliptical disc with $b = 0.5$. }
  \label{table:2}
\end{table}

\vspace{-1.em}

\begin{table}[H]
  \centering
  \begin{tabular}[H]{ccccc}
 & $I_{\sigma_0}$&  $I_{\sigma_1}$&  $I_{\sigma_2}$&  $I_{\sigma}$  \\
\hline\\[-0.75em]
$\varepsilon_{\textrm{rel}}$ & $0.055 \promille$ & $0.195\promille$ & $9.1 \times 10^{-4} \promille $ & $0.334\promille$ 
  \end{tabular}
  \caption{\small Maximum value of the relative error $\varepsilon_{\textrm{rel}} = \max \vert I^\textrm{exact}_{\sigma_i} - I^\textrm{comp}_{\sigma_i} \vert $ in per mil for an elliptical disc for $b = 0.5$ and $a=0.5:0.05:1.5$. }
  \label{table:3}
\end{table}

\section*{Acknowledgements}
The authors would like to express their thanks to  A. Bendali (INSA) for fruitful discussions and M. Fares (CERFACS) for the numerical computations achieved with the CERFACS code CESC. Part of this work was supported by the French National Research Agency under grant no. ANR$-08-$SYSC$-001$.

\appendix
\section{The case of an unit circle disc}
\label{app}
Let $\sC$ be the circle with radius $1$. We aim here at computing the integral
\begin{equation*}
   I_{\sC}\;=\;\dfrac{1}{4\pi}\int_{\sC}\dfrac{x_1\;y_1 }{\vert \bx-\by \vert} \: ds_\bx ds_\by
\end{equation*}
   This integral is rewritten in polar coordinates ($r$, $\phi$ for $\bx$ and $\rho$, $\phi'$ for $\by$) as
   \begin{equation*}
   I_\sC\;=\;\frac{1}{4\pi}\int_0^1\int_0^{2\pi}\int_0^1\int_0^{2\pi}\dfrac{r\cos\phi\;\rho\cos\phi' }{\sqrt{r^2+\rho^2-2r\rho\cos(\phi-\phi')}} \: rdrd\phi \;\rho d\rho d\phi'
   \end{equation*}
and evaluated using to the formula $(3.4.5)$ of \cite[p70]{copson1947}
\begin{equation*}
\! \int_0^{2 \pi} \dfrac{\cos \phi \:  d\phi}{\sqrt{\rho^2+r^2-2r\rho \cos(\phi-\phi')}} = \dfrac{4\cos\phi' }{\rho r} \!  \int_0^{ \min(\rho,r)}\dfrac{t^{2} \: dt}{\sqrt{\rho^2-t^2}\sqrt{r^2-t^2}}.
\end{equation*}
This leads to
   \begin{equation*}
   I_\sC\;=\;\frac{1}{\pi}\int_0^1\int_0^1\int_0^{ \min(\rho,r)} \dfrac{r\rho \: t^{2} }{\sqrt{\rho^2-t^2}\sqrt{r^2-t^2}}\: dr d\rho dt \ds\int_0^{2\pi}\cos^2\phi'd\phi'.
   \end{equation*}
 This integral is symmetric in $\rho$ and $r$, and since $\ds \int_0^{2\pi} \cos^2 \phi' d \phi' = \pi$, we have
\begin{equation*}
   I_{\sC}  =  2  \int_{t=0}^{1} t^2 \int_{\rho=t}^{1} \dfrac{\rho }{\sqrt{\rho^2-t^2}}\int_{r=\rho}^{1}   \dfrac{r }{\sqrt{r^2-t^2}}\: dr \: d\rho \: dt=\frac{2}{15}.
\end{equation*}





\bibliographystyle{elsarticle-num}
\bibliography{biblio-acoustique}

\begin{thebibliography}{1}
\expandafter\ifx\csname url\endcsname\relax
  \def\url#1{\texttt{#1}}\fi
\expandafter\ifx\csname urlprefix\endcsname\relax\def\urlprefix{URL }\fi
\expandafter\ifx\csname href\endcsname\relax
  \def\href#1#2{#2} \def\path#1{#1}\fi

\bibitem{leppington1973}
F.~Leppington, H.~Levine, {Reflexion and transmission at a plane screen with
  periodically arranged circular or elliptical apertures}, J. Fluid Mech 61
  (1973) 109--127.

\bibitem{laurens2012}
S.~Laurens, S.~Tordeux, A.~Bendali, M.~Fares, R.~Kotiuga, {Lower and upper
  bounds for the Rayleigh conductivity of a perforated plate}, submitted,
  http://hal.archives-ouvertes.fr/hal-00686438.

\bibitem{copson1947}
E.~Copson, {On the problem of the electrified disc}, Proceedings of the
  Edinburgh Mathematical Society (Series 2) 8~(01) (1947) 14--19.

\end{thebibliography}







\end{document}